\documentclass{commat}

\usepackage{stackengine}
\stackMath
 
\savestack{\rrl}{$\mathbin{\stackanchor[5pt]{\rightrightarrows}{\leftarrow}}$}

\title{%
    Pre-crossed modules and rack homology
    }

\author{%
    Jacob Mostovoy
    }

\affiliation{
    \address{Jacob Mostovoy --
    Departamento de Matem\'aticas, CINVESTAV-IPN, Av.\ Instituto Polit\'ecnico Nacional 2508, Col.\ San Pedro Zacatenco, Ciudad de M\'exico. 
        }
    \email{%
    jacob@math.cinvestav.mx
    }
    }

\abstract{%
    We define a homology theory for pre-crossed modules that specifies to rack homology in the case when the pre-crossed module is freely generated by a rack. 
    }



\VOLUME{30}
\NUMBER{2}
\YEAR{2022}
\firstpage{119}
\DOI{https://doi.org/10.46298/cm.10153}

\begin{paper}

\section{Introduction} 

The fundamental method of the combinatorial topology consists in approximating a topological space by a sequence of its $k$-dimensional skeleta. Depending on the chosen combinatorial model for the cells, these skeleta may be cubical or simplicial or of some other type and may satisfy different combinatorial properties but, invariably, the 0th approximation is a discrete set and the first approximation is a graph. If a topological space is equipped with a product (or some other algebraic structure) that is compatible with the skeleta, its restriction to the $k$-skeleton can give rise to a new algebraic structure. When $k=0$, this algebraic structure is just the discrete version of the structure on the space: for instance, the vertices of a simplicial group form a discrete group. When $k=1$, we get genuinely new algebraic structures: simplicial groups give rise to pre-crossed modules, pre-cubical products produce augmented racks and, as a linearized version of the same construction, DG Lie algebras give us (augmented) Leibniz algebras.

It has been observed (with some hard work involved) that racks and Leibniz algebras come with a new kind of homology. These invariants were first defined intrinsically, and then it was realized that they have a natural interpretation in terms of the corresponding cubical and DG structures. The goal of the present note is to point out that there is a similar homology for pre-crossed modules. Since both pre-crossed modules and augmented racks are one-dimensional ``shadows" of associative products, the former being simplicial and the latter pre-cubical, it is very natural that these new invariants are related to the rack homology.
 
The definition of the ``pre-crossed'' homology is similar in spirit to the definition of the Leibniz homology via DG Lie algebras: it also uses the left adjoint to the 1-truncation functor, see \cite{MDGLie}. In fact, Loday, in his search of the integration procedure for Leibniz algebras, wrote that  ``pre-crossed modules of groups give rise to a notion which almost fulfill the requirements about Leibniz algebras'' \cite[page 124]{LodayK}.  However, 
the ``pre-crossed'' homology of an abelian group is different from what Loday expected from his conjectural Leibniz homology for groups and we do not explore the possible connections with the Leibniz algebras here. 
  
Recall that an {augmented rack} with the {augmentation group} $G$ is a right $G$-set $X$ together with a morphism of $G$-sets $$X\xrightarrow{\pi} G$$ where $G$ acts on itself by conjugation.  If $X$ is a group on which $G$ acts by automorphisms and $\pi$ is a homomorphism, the augmented rack $X\xrightarrow{\pi} G$ is called a {pre-crossed module}. We will write the action of $G$ on $X$  as $$x\cdot g= x^g.$$
An augmented rack $X\xrightarrow{\pi} G$ gives rise to the pre-crossed module $$F(X)\xrightarrow{\overline{\pi}} G,$$ where $F(X)$ is the free group on $X$ and $\overline{\pi}$ is induced by $\pi$. We will show (Theorem~\ref{ra}) that the ``pre-crossed'' homology of this pre-crossed module  coincides with the rack homology of $X\xrightarrow{\pi} G$. 
Another case when we can identify the ``pre-crossed'' homology is that of a pre-crossed module $X\xrightarrow{\pi} G$ with a trivial action: it turns out to be the tensor algebra on the homology of the group $X$. In other cases,  these new invariants remain elusive, not unlike the rack homology. 

The tools that we use are classical. In order to relate the the rack homology to the ``pre-crossed'' homology we use the group-completion theorem of Quillen \cite{Q}. As for the computation of the homology for the pre-crossed modules with the trivial action (Theorem~\ref{tri}), it rests on the identification of the classifying space for the ``pre-crossed'' homology as a certain twisted version of the Milnor-Carlsson construction of a circle \cite{Carl}. When the action is trivial, the twisting vanishes and we use the result of Carlsson which describes the homotopy type of this space.
 
We provide no background on racks and their homology since extensive literature exists on the subject; see, for instance \cite{FRS}. For details on simplicial groups we refer to \cite{Wu}.

\section{The pre-crossed homology}

\subsection{The universal simplicial envelope of a pre-crossed module}\label{use}
Pre-crossed modules have been studied in some detail in the context of the truncation functors for simplicial groups. 

A  \emph{simplicial group truncated to order $p$} is a finite sequence of groups $G_p, \ldots, G_0$ together with the face and degeneracy maps between them that satisfy the usual simplicial identities; see \cite{Cond}. 
Simplicial groups truncated to order $1$ can be identified with pre-crossed modules. 
Indeed, let $G_1\ \rrl\ G_0$ be a simplicial group truncated to order $1$ with face maps $d_0, d_1$ and the degeneracy $s_0$.  The group $G_0$ can be identified with the subgroup $s_0G_0$ of the degenerate simplices in $G_1$ and, therefore, acts on $G_1$ by conjugation. Since $\ker d_0\subseteq G_1$ is a normal subgroup, $G_0$ acts on it by automorphisms. Therefore, $d_1: \ker{d_0}\to G_0$ is a pre-crossed module.

Conversely, given a pre-crossed module $X\xrightarrow{\pi} G$, one has the simplicial group truncated to order $1$  
$$X\rtimes G\ \rrl\ G$$ with 
\begin{align*}
d_0(x,g)&=g,\\ d_1(x,g)&=\pi(x) g
\end{align*} and 
$$s_0(g)=(1,g).$$
We see that the categories of pre-crossed modules and simplicial groups truncated to order $1$ are equivalent. In what follows we will identify these categories.

The \emph{truncation functor} from the category of simplicial groups to the category of simplicial groups truncated to order $p$ consists in discarding the simplices in dimensions higher than $p$. The truncation functor has the right and the left adjoint functors, known as the \emph{$p$-coskeleton} and the \emph{$p$-skeleton} functors. 
There is a canonical map from the $p$-skeleton to the $p$-coskeleton whose existence is guaranteed by the fact that they are the two adjoints to the truncation functor.

The 1-skeleton $E_*(X,G,\pi)$ of a pre-crossed module $X\xrightarrow{\pi} G$ is constructed as follows. 

Let $Y_*$ be the simplicial set which in degrees lower than 2 coincides with $X\rtimes G\ \rrl\ G$, and in degrees 2 and higher is obtained by freely adding degenerate simplices. Define the simplicial group $E_*(X,G,\pi)$ to be the quotient of the free simplicial group $F(Y_*)$ generated by $Y_*$ by the relations of the form
$$a*b = a\cdot b,$$
for each $a,b\in Y_1$, where the product on the left-hand side is taken in $F(Y_*)$ and on the right-hand side in $Y_1$. We call this group the \emph{universal simplicial envelope} of the pre-crossed module $\pi$. 
The simplicial group $E_*(X,G,\pi)$ in degree $k$ is generated by 
$$(x,g)_j = s_0^{j}s_1^{k-1-j}(x,g),$$ 
where $(x,g)\in Y_1$ and $j=0,\ldots, k-1$, and by $$g \in s_0^k G.$$ 

With these generators, the relations can be written as
$$(x,g)_j =   (x, 1)_j *g = g* (x^{g}, 1)_j, $$
$$(x_1, 1)_j* (x_2,1)_j = (x_1 x_2,1)_j,$$
$$g*h=gh,$$
and $$(1,1)_i=(1,1)_j=1.$$
It follows that $E_k(X,G,\pi)$ is a semi-direct product of the free product of $k$ copies of $X$ with a copy of the group $G$ (see Remark~2.5 in \cite{Cond}).

The faces and degeneracies in $E_*(X,G,\pi)$ 
restrict to the identity map on the semi-direct factor $G$. 
For $(x,1)_j\in E_k(X,G,\pi)$ we have 
\begin{equation}\label{e1}
s_i (x,1)_j = \left\{ \begin{array}{ll}(x,1)_{j+1}&\ \text{when}\ i\leq j\\(x,1)_{j}&\ \text{when}\ i>j\end{array} \right.
\end{equation}
and

\begin{equation}\label{e2}
d_i (x,1)_j = \left\{
\begin{array}{ll}
\quad 1& \ \text{when}\ i=j=0,\\ 
(x,1)_{j}&\ \text{when}\ i>j\ \text{and}\ j<k-1\ ; \\
(x,1)_{j-1}&\ \text{when}\ i\leq j\ \text{and}\ j>0; \\
\ \pi(x)& \ \text{when}\ j=k-1,\ i=k.\end{array} \right.
\end{equation}

\subsection{The pre-crossed homology}

From the explicit constructions in the previous subsection, we see that the group $G$ of the pre-crossed module $X\xrightarrow{\pi} G$  acts on 
 $E_k(X,G,\pi)$ freely on the right by multiplication. 
Write $$PH_m(X,G,\pi) := H_m(E_*(X,G,\pi)/ G).$$
We call $PH_m(X,G,\pi)$ the $m$th 
\emph{pre-crossed} homology group of the pre-crossed module $X\xrightarrow{\pi} G$.
 
One can exclude any mention of the group $G$ in the definition of the pre-crossed homology.    
Define a \emph{pre-crossed action of $X$ on itself} as a homomorphism $$\phi: X\to \mathrm{Aut}(X)$$ such that $\phi: X\to \phi(X)$ is a pre-crossed module. A pre-crossed module $X\xrightarrow{\pi} G$ gives rise to a pre-crossed action $\phi_\pi$ of $X$ on itself, namely, the composition
$$X\xrightarrow{\pi} G\to \mathrm{Aut}(X).$$ 
Indeed, $\pi$ send the action of $G$ on $X$ to conjugation on $G$ and any surjective homomorphism sends conjugation to conjugation.

If $\phi=\phi_\pi$, by the definition of the universal simplicial envelope we have that
$$E_*(X,G,\pi)/ G = E_*(X,\phi(X),\phi)/ \phi(X),$$
so that the pre-crossed homology of $X\xrightarrow{\pi} G$ and $X\xrightarrow{\phi} \phi(X)$ are the same. In this manner, the pre-crossed homology does not depend on $G$ but only on the group $X$ and the homomorphism $\phi: X\to \mathrm{Aut}(X)$. 

\begin{remark} A pre-crossed action of a group on itself is a natural ``non-augmented'' version of a pre-crossed module. Essentially, it is the same thing as a rack with a compatible group structure.
\end{remark}

Pre-crossed modules form a category and the pre-crossed homology is easily seen to be functorial.
 
\begin{proposition} There is a natural transformation of functors  $$PH_m(X, G,\pi)\to H_m(G)$$
from the pre-crossed homology of $\pi:X\to G$ to the usual homology of the group $G$.
\end{proposition}

\begin{proof}
Write $M_*(X,G,\pi)$ for  the 1-coskeleton of the simplicial group truncated to order one $X\rtimes G\ \rrl\ G$. Since $G$ appears as a subgroup of the degenerate simplices in $M_k(X,G,\pi)$ for all $k$, it acts on $M_*(X,G,\pi)$ freely on the right. The canonical map of the 1-skeleton to the 1-coskeleton gives rise to a functorial map
$$H_m(E_*(X,G,\pi)/ G)\to H_m(M_*(X,G,\pi)/ G).$$
Without loss of generality, we can assume that $\pi$ is surjective. It follows from \cite[Theorem~1.3]{Cond}, that the Moore complex of $M_*(X,G,\pi)$ is acyclic when $\pi$ is surjective and, therefore, $M_*(X,G,\pi)$ is contractible.
Therefore $H_m(M_*(X,G,\pi)/ G)$ is the usual homology of $G$.
\end{proof}

\begin{remark} 
As we have already mentioned, the relationship between the pre-crossed and usual homology mirrors the relationship between the Leibniz and the Lie algebra homology in the context of DG Lie algebras, see \cite{MDGLie}. There are, in fact, other cases when a truncation (forgetful) functor has both right and left adjoints which lead to two kinds of (co)homologies with a natural transformation between them. The oldest such example is, probably, that of Lie triple systems that have two different enveloping Lie algebras, the universal and the standard enveloping algebra, which give rise to two versions of cohomology for Lie triple systems, see \cite{Ha},  \cite{Ho}, \cite{Li}. In fact, the Yamaguti's cohomology of Lie triple systems \cite{Y} may have been the first (covert) appearance of Leibniz cohomology on the mathematical scene, see \cite{Z}.
\end{remark}

There are two  cases of the pre-crossed homology that may be of particular interest: the homology of the pre-crossed module $F(X)\to G$ associated with the augmented rack $X\to G$ and the homology of the pre-crossed module $\mathrm{id}: G\to G$, where $G$ is a group acting on itself by conjugation. We will see in the next section that the first of these homologies coincides with the usual rack homology. Computing the latter is outside the scope of this note, although we explicitly identify it in the case of abelian $G$.

\begin{remark} 
There are many functorially defined simplicial groups that fit in between the 1-skeleton and the 1-coskeleton of a pre-crossed module. For instance, if $R_*(0)$ is the kernel of the natural homomorphism
$$E_*=E_*(X,G,\pi)\to M_*(X,G,\pi),$$
define for all integer $n>0$ $$R_*(n)=[E_*, R_*(n-1)].$$
Then $$E_*(n)=E_*/R_*(n)$$ is a simplicial group whose 1-truncation gives the pre-crossed module $X\xrightarrow{\pi} G$ and which carries an action of $G$. One may then consider the homology of $E_*(n)/G$ as a functor that interpolates between the usual homology of $G$ and the pre-crossed homology of $X\xrightarrow{\pi} G$. Similar functors exist for Leibniz algebras \cite{MDGLie}.
 \end{remark}

\section{Relationship with the rack homology}

\subsection{The Clauwens monoid and the rack space}

Each augmented rack $X\xrightarrow{\pi} G$ embeds into a topological monoid that we call the \emph{Clauwens monoid} of $X\xrightarrow{\pi} G$ (see \cite{Clau}). One way to describe this construction is as follows (see~\cite{MRacks}).

Consider the one-dimensional cell complex $\Gamma_\pi$ with the set of vertices $G$ and the set of edges $X\times G$: 
for each $x\in X$ with $\pi(x)= h g^{-1}$ the complex $\Gamma_\pi$ has one edge $(x,g)$ from $g$ to $h$. The two-sided action of $G$ on itself by multiplications extends to the action of $G$ on $\Gamma_\pi$:
$$\quad h\cdot (x, g) = (x^{h^{-1}}, hg)$$
and
$$(x, g) \cdot h = (x, gh).$$
The \emph{Clauwens monoid} of $X\xrightarrow{\pi} G$ is obtained from $\Gamma_\pi$ by taking the quotient of the free monoid (that is, the James reduced product) on $\Gamma_\pi $ by the relations 
$$g*t = g\cdot t\quad \text{and}\quad t*g = t\cdot g,$$ 
where $g\in G\subset  \Gamma_\pi, t\in \Gamma_\pi$, the star denotes the product in the free monoid and the dot -- the two-sided action of $G$ on $\Gamma_\pi$. 

There is a free right action of $G$ on the Clauwens monoid and the quotient is known as the \emph{rack space} of $X\xrightarrow{\pi} G$. The  \emph{rack homology} of $\pi$ is then defined as the homology of the rack space of $X\xrightarrow{\pi} G$. 

\subsection{The group completion of the Clauwens monoid}

The Clauwens monoid and the rack space have a natural structure of cubical complexes. They can, in fact, also be described in terms of simplicial groups.

\begin{theorem}\label{ra}
Let $X\xrightarrow{\pi} G$ be an augmented rack and $F(X)\xrightarrow{\overline{\pi}} G$ be its associated pre-crossed module. Then, the geometric realization of the simplicial group $E_*(F(X), G, \overline{\pi})$ is homotopy equivalent to the Clauwens monoid of $X\xrightarrow{\pi} G$ and the geometric realization of $E_*(F(X), G, \overline{\pi})/G$ is the rack space of $\pi$.
\end{theorem}

The proof is, essentially, an application of Quillen's Theorem Q.4 in \cite{Q} on the group completion of a simplicial monoid. This theorem implies that if $M$ is a connected simplicial monoid and the canonical homomorphism of $M$ to its group completion $\overline{M}$ induces 
a homotopy equivalence of the classifying spaces $BM\to B\overline{M}$, the completion map $M\to \overline{M}$ itself is a homotopy equivalence\footnote{Theorem~Q.4 gives an isomorphism in homology, rather than homotopy, but this is sufficient since a connected monoid is homotopy simple.} after geometric realization. In order to use this result, we need the following:

\begin{lemma}\label{above}
Assume that a group $G$ acts on a free monoid $F$ by automorphisms. Then, the natural map 
$$B(F\rtimes G)\to  B(\overline{F\rtimes G})$$ is a homotopy equivalence. 
\end{lemma}

\begin{proof} 
By Proposition~Q.1 of \cite{Q}, $BF\to  B\overline{F}$ is a homotopy equivalence. Also, we have 
$$\overline{F\rtimes G}=  \overline{F}\rtimes G.$$
On the other hand, $B(F\rtimes G)$ and $B(\overline{F}\rtimes G)$ are fibre bundles over $BG$ and the natural map 
$$B(F\rtimes G) \to B(\overline{F}\rtimes G)$$ is a homotopy equivalence on the fibers and, therefore, a homotopy equivalence.
\end{proof}

We will also need Proposition~Q.2 of \cite{Q} which says that if, for a simplicial monoid $M$ each $M_k$ has the property that $BM_k\to B\overline{M}_k$ is a homotopy equivalence, the same holds true for $M$ itself.

\begin{proof}[Proof of Theorem~\ref{ra}]
Let $Y_*$ be the simplicial set which in degrees lower than 2 coincides with $X\rtimes G\ \rrl\ G$, and in degrees 2 and higher is obtained by freely adding degenerate simplices; it has been introduced in Section~\ref{use}. (It has been defined in the situation when $X$ is a group but its definition does not use the group structure on $X$). The geometric realization of $Y_*$ is precisely the cell complex $\Gamma_\pi$.

Define the simplicial Clauwens monoid $C_*$ as the quotient of the free monoid on $Y_*$ by the relations 
$$g*h=gh, \ \ g*(x,h) = (x^{g^{-1}}, gh) \text{\ and\ } (x,g)*h = (x, gh),$$
for $g,h\in G, x\in X$.
 
The geometric realization of the simplicial Clauwens monoid is homeomorphic to the topological Clauwens monoid since the operations of forming a free monoid and applying the relations can be expressed in terms of colimits and finite limits and, therefore, commute with the geometric realization. 

The group completion of $C_*$ coincides with $E_*(F(X),G,\pi)$. Each $C_k$ is a semi-direct product of $G$ and a free monoid so, by Lemma~\ref{above}, Proposition~Q.1 and Theorem~Q.4 of \cite{Q}, the geometric realization of $C_*$ is homotopy equivalent to that of $E_*(F(X),G,\pi)$.

Since the action of $G$ on $C_*$ and on $E_*(X,G,\pi)$ is free, the rack space of $X\xrightarrow{\pi} G$, which is the realization of $C_*/G$, is also homotopy equivalent to that of $E_*(X,G,\pi)/G$.
\end{proof}

\section{Pre-crossed modules with the trivial action and the Milnor-Carlsson construction}

When we compute $PH_m(X,G,\pi)$ with the trivial action of $G$, the group $G$ can be replaced by the trivial group. In this case, $PH_m(X,G,\pi)$ coincides with the homology of the universal simplicial envelope, which can be explicitly identified: in degree $k$, it is the free product of $k$ copies of $X$, with the faces and degeneracies as in (\ref{e1}),  (\ref{e2}).  

This simplicial group is well-known: it coincides with the based Milnor-Carlsson construction $J_X(S^1,*)$ on the circle $S^1$ with the coefficients in $X$ and the trivial action, as defined in \cite{Carl}. Recall that the original Milnor construction is the simplicial group version of the James reduced product: it assigns to a simplicial set $S$ the free simplicial group generated by $S$. Its homotopy type is that of the loop space on the suspension of $S$. Carlsson generalized it in \cite{Carl} to simplicial sets with a group action. Namely, if $S$ is a simplicial set with an action of a group $X$, 
Carlsson defines $J_X(S,*)$ as the simplicial group whose group of $n$-simplices 
has generators $[s,x]$ with $s\in S_n$, the set of $n$-simplices of $S$  and $x\in X$, modulo the relations
$$[s, x][s^x, y]=[s, xy]$$
and
$$[*, x]=1,$$
where $*$ is the $n$-fold degeneracy of the base-point. The faces and the degeneracies of  $J_X(S,*)$ are induced by those of $S$.

Let $S^1$ be the simplicial circle with one non-degenerate simplex $c_0$ in degree 0 and $c_1$ in degree 1; for $n>1$, the $n$-simplices of $S^1$ are all degenerate and form the $(n+1)$-element set
$$\{ *=s_0^n c_0,  s_0^{n-1} c_1, s_0^{n-2}s_1 c_1, \ldots, s_1^{n-1} c_1\}.$$
We are interested in the case $S=S^1$ with the trivial action of $X$ and $G=1$, the trivial group. In this case, $J_X(S^1,*)$ is identified with $E_*(X,1,\pi)$ by
$$[  s_0^{j} s_1^{n-1-j} c_1, x] \mapsto (x,1)_j.$$ 

The homotopy type of $J_X(S,*)$ has been determined by Carlsson in \cite{Carl}. In particular, for a simplicial set $S$ with the trivial action of $X$, the simplicial set $J_X(S,*)$ has the homotopy type of the loop space on the smash product $S\wedge BX$ of $S$ with the classifying space of $X$. In our case, 
$$J_X(S^1,*)\simeq \Omega\Sigma BX,$$ 
the loop space on the suspension of the classifying space of $X$.
By the Bott-Samelson Theorem \cite{BS}, its homology with coefficients in a field is the tensor algebra on the homology of $X$ with the same coefficients.

As a corollary, we have
\begin{theorem}\label{tri} When the action of $G$ on $X$ is trivial,
$$PH_*(X,G,\pi)\otimes \mathbb{Q} = T_*(H_*(X,\mathbb{Q})).$$
\end{theorem}
 
When the action in a pre-crossed module is non-trivial, the universal simplicial envelope may be thought of as a ``twisted'' version of $J_X(S^1,*)$; hopefully, this point of view may lead to the identification of the pre-crossed homology in more interesting cases.


\EditInfo{October 15, 2022}{December 27, 2022}{Ivan Kaygorodov}

\end{paper}